\documentclass[12pt,reqno]{article}

\usepackage[usenames]{color}
\usepackage{amssymb}
\usepackage{amsmath}
\usepackage{amsthm}
\usepackage{amsfonts}
\usepackage{amscd}
\usepackage{graphicx}

\usepackage[colorlinks=true,
linkcolor=webgreen,
filecolor=webbrown,
citecolor=webgreen]{hyperref}

\definecolor{webgreen}{rgb}{0,.5,0}
\definecolor{webbrown}{rgb}{.6,0,0}

\usepackage{color}

\usepackage{graphics}
\usepackage{latexsym}

\begin{document}

\theoremstyle{plain}
\newtheorem{theorem}{Theorem}
\newtheorem{corollary}[theorem]{Corollary}
\newtheorem{lemma}[theorem]{Lemma}
\newtheorem{proposition}[theorem]{Proposition}

\theoremstyle{definition}
\newtheorem{definition}[theorem]{Definition}
\newtheorem{example}[theorem]{Example}
\newtheorem{conjecture}[theorem]{Conjecture}

\theoremstyle{remark}
\newtheorem{remark}[theorem]{Remark}

\begin{center}
\vskip 1cm{\LARGE\bf On a New Alternating Convolution Formula for the Super Catalan Numbers
\vskip 1cm}
\large
Jovan Miki\'{c}\\
University of Banja Luka\\
Faculty of Technology\\
Republic of Srpska\\
\href{mailto:jovan.mikic@tf.unibl.org}{\tt jovan.mikic@tf.unibl.org} \\
\end{center}

\vskip .2in

\begin{abstract}
We present a new alternating convolution formula  for the super Catalan numbers  which arises as a  generalization of two known binomial identities. We prove a generalization of this formula  by using auxiliary sums, recurrence relations, and induction.  By using a new method,  we prove one interesting divisibility result with super Catalan numbers.

\end{abstract}

\section{Introduction}

Let $l$ be a fixed non-negative integer, and let $n$ be an arbitrary non-negative integer.  Let $C_n=\frac{1}{n+1}\binom{2n}{n}$ denote the $n$th Catalan number.

Catalan ($1874$) observed \cite {EC} that numbers $S(n,l)=\frac{\binom{2n}{n}\binom{2l}{l}}{\binom{n+l}{n}}$  are always integers.   See also \cite[p.\ 18]{Koshy}. In particular, $S(n,0)=\binom{2n}{n}$ and $S(n,1)=2C_n$.  Gessel \cite [p.\ 11]{IG} referred to these numbers as the super Catalan numbers, since $\frac{1}{2}S(n,1)$ is the Catalan number $C_n$.

We shall call $S(n,l)$ as the $n$th super Catalan number of order $l$.  By symmetry, $S(n,l)=S(l,n)$. Also $S(n,l)$ is equal to $\frac{(2n)!(2l)!}{n!\cdot l!\cdot (n+l)!}$.

It is known that  $S(n,l)$ is always an even integer except for the case $n=l=0$. See \cite[Introduction]{EAllen} and \cite[Eq.~(1), p.\ 1]{DC}.  More generally, for positive $l$, numbers  $\frac{1}{2} S(n,l)$ can be viewed as special cases of super ballot numbers \cite[p.\ 11]{IG}.

 For only a few values of $l$, there exist combinatorial interpretations of $S(n,l)$.  See, for example, \cite{EAllen, DC, ChenWang, gesselxin04, pippenger, schaeffer03}. The problem of finding a combinatorial interpretation for super Catalan numbers of an arbitrary order $l$ is an intriguing open problem.

There are several binomial coefficient identities for super Catalan numbers.  For example,  the identity of von Szily ($1894$): \cite[Eq.~(29), p.\ 11]{IG}
\begin{equation}
S(n,l)=\sum_{k \in \mathbb{Z}}(-1)^k\binom{2n}{n+k}\binom{2l}{l+k}.\notag
\end{equation}

Note that  the identity of von Szily gives another proof  that the number $S(n,l)$ is always an integer. See also \cite[Eq.~(31); Eq.~(32), p.\ 12]{IG}.

We present the following alternating convolution formula for the super Catalan numbers:
\begin{theorem}\label{t:1}
For all non-negative integers $n$ and $l$, we have:
\begin{equation}
\sum_{k=0}^{2n}(-1)^k\binom{2n}{k}S(k,l)S(2n-k,l)=S(n,l)S(n+l,n).\label{eq:1}
\end{equation}
\end{theorem}

For $l=0$, Eq.~(\ref{eq:1}) reduces to a known binomial identity \cite[Example 3.6.2, p.\ 45]{MPWZ}:

\begin{equation}
\sum_{k=0}^{2n}(-1)^k{\binom{2n}{k}}\binom{2k}{k}\binom{4n-2k}{2n-k}={\binom{2n}{n}}^2.\label{eq:2}
\end{equation}

Recently, the  following  binomial coefficient identity involving
the Catalan numbers  \cite{JM3} was discovered:
\begin{equation}
\sum_{k=0}^{2n}(-1)^k\binom{2n}{k} C_k\, C_{2n-k}=C_n\binom{2n}{n}.\label{eq:3}
\end{equation}

See also \cite[Remark 23, p.\ 15]{JM2} and  \cite{WChu, prodinger}. For $l=1$, Eq.~(\ref{eq:1}) reduces to an identity equivalent to Eq.~(\ref{eq:3}). Therefore, we can see Eq.~(\ref{eq:1}) as a natural generalization of Eqns. (\ref{eq:2}) and (\ref{eq:3}).

Furthermore,  Eq.~(\ref{eq:1}) has the following generalization:

\begin{theorem}\label{t:2}
Let $n$, $l$, and  $t$ be non-negative integers such that $t\leq n $. Then 
\begin{equation}\label{eq:7}
\sum_{k=t}^{2n-t}(-1)^k\binom{2n-2t}{k-t}S(k,l)S(2n-k,l)=(-1)^t\frac{\binom{2l}{l}\binom{2t}{t}\binom{2(n+l-t)}{n+l-t}\binom{2n}{n}\binom{2n-2t}{n-t}}{\binom{n+l}{n}\binom{2n+l-t}{n}\binom{n}{t}}.
\end{equation}
\end{theorem}

When $t=0$, Eq.~(\ref{eq:7}) reduces to  Eq.~(\ref{eq:1}).

For $l=0$, Theorem \ref{t:2} reduces to a recently  \cite[Lemma 10, p.\ 7]{JM2} discovered formula:
\begin{equation}
\sum_{k=t}^{2n-t}(-1)^k\binom{2n-2t}{k-t}\binom{2k}{k}\binom{4n-2k}{2n-k}=(-1)^t\frac{\binom{2n}{n}\binom{2t}{t}\binom{2n-2t}{n-t}}{\binom{2n-t}{t}}.\label{eq:8}
\end{equation}

It is readily verified that for $l=1$, Theorem \ref{t:2} is equivalent to another recently   \cite[Lemma 13, p.\ 7]{JM2} discovered formula:
\begin{equation}\label{eq:9}
\sum_{k=t}^{2n-t}(-1)^k\binom{2n-2t}{k-t}C_kC_{2n-k}=(-1)^t\frac{C_n\binom{2t}{t}\binom{2n-2t}{n-t}}{\binom{2n+1-t}{t}}.
\end{equation}

We prove Theorem \ref{t:2} by using auxiliary sums, recurrence relations, induction, and Eq.~(\ref{eq:8}).

Let us consider the following sum:
\begin{equation}
\varPsi(n,m,l)=\sum_{k=0}^{n}(-1)^k\binom{n}{k}^m S(k,l)S(n-k,l).\label{eq:10}
\end{equation}

Obviously,  $\varPsi(2n-1,m,l)=0$.

By Eq. (\ref{eq:1}), it follows that $\varPsi(2n,1,l)$ is divisible by $S(n,l)$. We establish the following theorem:

\begin{theorem}\label{t:3}
The sum $\varPsi(2n,m,l)$ is divisible by $S(n,l)$ for all non-negative integers $n$ and $l$ and  for all positive integers $m$.
\end{theorem}

For proving Theorem \ref{t:3}, we use Theorem \ref{t:2} and a method we call the ``method
of $D$ sums''.
\begin{definition}\label{def:1}
Let $n$, $j$, $l$, and $t$ be non-negative integers such that $j\leq \lfloor \frac{n}{2}\rfloor$, and let $m$ be a positive integer.
Let $A(n,m,l)=\sum_{k=0}^{n}\binom{n}{k}^m F(n,k,l)$,
where $F(n,k,l)$ is an integer-valued function that depends  on $n$, $k$, and $l$ (not on $m$). Then the $D$ sums for $A(n,m,l)$ are
\begin{equation}
D_A(n,j,t;l)=\sum_{u=0}^{n-2j}\binom{n-j}{u}\binom{n-j}{j+u}\binom{n}{j+u}^tF(n,j+u,l).\label{eq:11}
\end{equation}
\end{definition}
First, note that sum (\ref{eq:10}) is an instance of $A(n,m,l)$.
For $m\geq 2$, by Eq.~(\ref{eq:11}), it follows that
\begin{equation}\label{eq:12}
A(n,m,l)=D_A(n,0,m-2;l)\text{.}
\end{equation}
Furthermore, $D$ sums satisfy the following two recurrence relations \cite[Thm.\ 2, Thm.\ 3, p.\ 2]{JM1}:
\begin{align}
 D_A(n,j,t+1;l)&=\sum_{u=0}^{\lfloor \frac{n-2j}{2} \rfloor} \binom{n}{j+u}\binom{n-j}{u}D_A(n,j+u,t;l),\label{eq:13}\\
 D_A(n,j,0;l)&=\sum_{u=0}^{\lfloor \frac{n-2j}{2} \rfloor}\binom{n-j}{j+u}\binom{n-2j-u}{u}\sum_{v=0}^{n-2j-2u}\binom{n-2j-2u}{v}F(n,j+u+v,l).\label{eq:14}
\end{align}

\begin{definition}\label{def:2}
Let $F$ be from Definition \ref{def:1}, and let $n$, $l$, $t$ be non-negative integers such that $t\leq \lfloor \frac{n}{2}\rfloor$.
Then $A_t(n,l)$ denotes
\begin{equation}\label{eq:15}
\sum_{k=t}^{n-t}\binom{n-2t}{k-t}F(n,k,l)\mbox{.}
\end{equation}
\end{definition}

By substitution $k=u+j+v$, the inner sum of the right-side of (\ref{eq:14}) becomes
\begin{equation}\label{eq:16}
\sum_{v=0}^{n-2j-2u}\binom{n-2j-2u}{v}F(n,j+u+v,l)=A_{j+u}(n,l).
\end{equation}
It is readily verified \cite[Eq.~(1.4), p.~5]{Koshy} that $\binom{n-j}{j+u}\binom{n-2j-u}{u}=\binom{n-j}{u}\binom{n-j-u}{j+u}$. By using this fact and  Eq.~(\ref{eq:16}), Relation \ref{eq:14} becomes
\begin{equation}\label{eq:17}
 D_A(n,j,0;l)=\sum_{u=0}^{\lfloor \frac{n-2j}{2}\rfloor}\binom{n-j}{u}\binom{n-j-u}{j+u}A_{j+u}(n,l).
\end{equation}
From now on, for calculating $D_A(n,j,0;l)$ sum, we use Eq.~(\ref{eq:17}) instead of Relation \ref{eq:14}.

It turns out that  $\varPsi_t(2n,l)$ is the left-side of Theorem \ref{t:2}. By Theorem \ref{t:2}, it follows that
\begin{equation}
 \varPsi_t(2n,l)=(-1)^t\frac{\binom{2l}{l}\binom{2t}{t}\binom{2(n+l-t)}{n+l-t}\binom{2n}{n}\binom{2n-2t}{n-t}}{\binom{n+l}{n}\binom{2n+l-t}{n}\binom{n}{t}}.\label{eq:18}
\end{equation}

\section{Motivation}\label{sec:2}

In $ 1998$, Calkin  proved that the alternating binomial sum $\sum_{k=0}^{2n}(-1)^k\binom{2n}{k}^m $ is divisible by $\binom{2n}{n} $ for all non-negative integers $n$ and all positive integers $m$.   In 2007, Guo, Jouhet, and Zeng proved, among other things, two generalizations of Calkin's result \cite[Thm.\ 1.2, Thm.\ 1.3, p.\ 2]{VG}. Also, they presented an interesting congruence involving super Catalan numbers  \cite[Corollary 4.9., p.\ 11]{VG}.

The first application of $D$ sums \cite[Section 8]{JM1} was  for proving Calkin's result \cite[Thm.\ 1]{NC}. Also, by using $D$ sums, it was proved \cite {JM1} that $\sum_{k=0}^{2n}\binom{2n}{k}^m |n-k|$ is divisible by $n\binom{2n}{n}$ for all non-negative integers $n$ and all positive integers $m$. 

Recently, by the same method, it was proved \cite[Thm.\ 1]{JM2}  that $\sum_{k=0}^{2n}(-1)^k\binom{2n}{k}^m\binom{2k}{k}\binom{4n-2k}{2n-k}$ is divisible by $\binom{2n}{n}$ for all non-negative integers $n$ and all positive integers $m$. This result confirms Theorem \ref{t:3}  for $l=0$.

 Furthermore, it was proved  \cite[Corollary 4]{JM2}  that   $\sum_{k=0}^{2n}(-1)^k\binom{2n}{k}^m C_kC_{2n-k}$ is divisible by $\binom{2n}{n}$ for all non-negative integers $n$ and all positive integers $m$. This result is sharper than Theorem \ref{t:3} for $l=1$.  By Theorem \ref{t:3} and Remark \ref{r:3}, it can be shown that $\sum_{k=0}^{2n}(-1)^k\binom{2n}{k}^m C_kC_{2n-k}$  is divisible by $C_n$.

In order to prove Theorem \ref{t:3}, we give a slight generalization of a notion $D$ sum in Definition \ref{def:1}. The function F, from Definition \ref{def:1}, depends on $n$, $k$, and on new parametar $l$. In all previous cases \cite{JM1,JM2}, the function F depends only on $n$ and $k$.

Recently, Chu \cite[Theorem 5]{WChu} presented, among other, formula similar to Theorem \ref{t:2}. It is readily verified that his formula implies Eq.~(\ref{eq:2}) for the case  $\lambda=\mu=t=0$ and $n:=2n$. 
Theorems   \ref{t:1} and \ref{t:2} are not consequences of this formula.

This paper consists of two parts.

 In the first part, we prove Theorem \ref{t:2} by using recurrence relations, induction, and Eq.~(\ref{eq:8}). This proof of Theorem \ref{t:2} is interesting itself for two reasons. Firstly, the choice of auxiliary sums. We use four auxiliary sums similarly as in the proof of Eq.~(\ref{eq:8}) \cite[Corollary 22, p.\ 8]{JM2}.
Secondly, we obtain a relation between $\varPsi_t(2n-2,l+1)$ and $\varPsi_t(2n,l)$ and then use induction on $l$. 

In the second part, we prove Theorem \ref{t:3} by using method of $D$ sums. 

Our proof  of Theorem \ref{t:3} consists of two parts. Firstly, we show that $\varPsi(2n,2,l)$ is divisible by $S(n,l)$ by calculating $D_\varPsi(2n,j,0;l)$ and by using Eq.~(\ref{eq:12}). 
Then we show that $D_\varPsi(2n,j,1;l)$ is divisible by $S(n,l)$ for all non-negative integers $j$ and $n$ such that $j\leq n$. This result is sufficient to prove Theorem \ref{t:2} for $m \geq 3$. Namely, then by Relation \ref{eq:13} and induction, it can be shown that all $D_\varPsi(2n,j,t;l)$ are divisible by $S(n,l)$ for all integers $t\geq1$. By Relation \ref{eq:12}, it follows that $\varPsi(2n,m,l)$ is divisible by $S(n,l)$ for all $m \geq 3$. See \cite[Section 5]{JM1}.
In order to calculate $D_\varPsi(2n,j,0;l)$, we use Eq.~(\ref{eq:18}) and Relation (\ref{eq:17}).

The rest of the paper is structured as follows.
In Section \ref{sec:3}, we present our main lemma \ref{l:1}, auxiliary sums, and recurrence relations between these sums. These recurrences are given by lemmas \ref{l:2}, \ref{l:3}, and \ref{l:4}. For auxiliary sums, see \cite{JM4,JM2}.
In Section \ref{sec:4}, we prove all lemmas from Section \ref{sec:3}.
In Section \ref{sec:5}, we give a proof of Theorem \ref{t:2} by using the main lemma \ref{l:1}, induction, and Eq. (\ref{eq:8}).
In Section \ref{sec:6}, we give a proof of Theorem \ref{t:1} by using Theorem \ref{t:2}.
In Section \ref{sec:7}, we give a proof of Theorem \ref{t:3} by using Theorem \ref{t:2} and a method of $D$ sums.

\section{Main Lemma, Auxiliary Sums, and Other Lemmas}\label{sec:3}

Let $n$, $t$, and  $l$ be non-negative integers such that $t\leq \lfloor\frac{n}{2}\rfloor$.

By Definition \ref{def:1} and Eq.~(\ref{eq:10}), we know that, for the sum $\varPsi(n,m,l)$,
\begin{equation}
F(n,k,l)=(-1)^k S(k,l)S(n-k,l)\text{.}\notag
\end{equation}

Therefore, by Eq.~(\ref{eq:10}), Definition \ref{def:1}, Definition \ref{def:2}, and Eq.~(\ref{eq:15}), it follows that 
\begin{equation}
 \varPsi_t(n,l)=\sum_{k=t}^{n-t}(-1)^k\binom{n-2t}{k-t}S(k,l)S(n-k,l)\label{eq:19}.
\end{equation}

By changing $k$ to $n-k$ in Eq.~(\ref{eq:19}), it follows that 
\begin{equation}
 \varPsi_t(2n-1,l)=0.\label{eq:20}
\end{equation}

For $l=0$, by Eq.~(\ref{eq:8}), it follows that
\begin{equation}
 \varPsi_t(2n,0)=(-1)^t\frac{\binom{2n}{n}\binom{2t}{t}\binom{2n-2t}{n-t}}{\binom{2n-t}{t}}.\label{eq:21}
\end{equation}

For $t=n$, by Eq.~(\ref{eq:19}), it follows that 
\begin{equation}
\varPsi_n(2n,l)=(-1)^nS(n,l)^2.\label{eq:22}
\end{equation}

Now we are ready for the main lemma.
\begin{lemma}\label{l:1}
 Let $n$ be a positive integer, and let $t$ and  $l$ be non-negative integers such that $t<n$.  Then the following recurrence relation is true:
\begin{equation}
 \varPsi_t(2n-2,l+1)=\frac{(2n+l-t)(2l+1)}{ 2(2n-1-2t)(2n-1)} \varPsi_t(2n,l).\label{eq:23}\\
\end{equation}
\end{lemma}

In order to prove lemma \ref{l:1}, we introduce four auxiliary sums, as follows:
\begin{align}
P_t(n,l)&=\sum_{k=t}^{n-t}(-1)^k(n-t-k)\binom{n-2t}{k-t}S(k,l)S(n-k,l)\label{eq:24},\\
R_t(n,l)&=\sum_{k=t}^{n-t}(-1)^k \frac{2l+1}{k+l+1}\binom{n-2t}{k-t}S(k,l)S(n-k,l)\label{eq:25},\\
R'_t(n,l)&=\sum_{k=t}^{n-t}(-1)^k \frac{2l+1}{(k+l+1)(n-k+l+1)}\binom{n-2t}{k-t}S(k,l)S(n-k,l)\label{eq:26},\\
T_t(n,l)&=\sum_{k=t}^{n-t}(-1)^k \,\frac{(n-t-k)(2l+1)}{k+l+1}\binom{n-2t}{k-t}S(k,l)S(n-k,l)\label{eq:27}.
\end{align}

It is readily verified that 
\begin{align}
P_t(2n,l)&=(n-t)\varPsi_t(2n,l),\label{eq:28}\\
R'_t(2n,l)&=\frac{1}{n+l+1}R_t(2n,l).\label{eq:29}
\end{align}

For example,  the proof of Eq.~(\ref{eq:28}) is as follows.

By Eq.~(\ref{eq:24}),
\begin{equation}
P_t(2n,l)=\sum_{k=t}^{2n-t}(-1)^k(2n-t-k)\binom{2n-2t}{k-t}S(k,l)S(2n-k,l).\notag
\end{equation}

By changing $k$ to $2n-k$ in the above equation and use symmetry of binomial coefficients, we obtain
\begin{equation}
P_t(2n,l)=\sum_{k=t}^{2n-t}(-1)^k(k-t)\binom{2n-2t}{k-t}S(k,l)S(2n-k,l).\notag
\end{equation}

If we add two above equations and use Eq.~(\ref{eq:19}),  Eq.~(\ref{eq:28})  follows. The proof of Eq.~(\ref{eq:29}) is similar to the proof of Eq.~(\ref{eq:28}).

We present other three lemmas.

\begin{lemma}\label{l:2}
Let $n$, $l$, and $t$ be non-negative integers.
Then the following recurrence relation is true:
\begin{equation}
R_t(2n,l)=\frac{n+l+1}{4(2l+1)}\,\varPsi_t(2n,l+1).\label{eq:30}
\end{equation}
\end{lemma}

\begin{lemma}\label{l:3}
Let $n$ be a positive integer, and let $l$ and $t$ be non-negative integers such that $t<n$.
Then the following recurrence relation is true:
\begin{equation}
\varPsi_t(2n,l)=4R_t(2n-1,l).\label{eq:31}
\end{equation}
\end{lemma}

\begin{lemma}\label{l:4}
Let $n$ be a positive integer, and let $l$ and $t$ be non-negative integers such that $t<n$.
Then the following recurrence relation is true:
\begin{equation}\label{eq:32}
(2n+l-t)R_t(2n-1,l)=\frac{2(2n-1-2t)(2n-1)}{n+l}R_t(2n-2,l).
\end{equation}
\end{lemma}

By lemmas \ref{l:2}, \ref{l:3}, and \ref{l:4}, the main lemma  \ref{l:1} follows. We shall use $T_t(n,l)$ only for the proof of lemma \ref{l:4}.

Note that, for $l=0$, Eqns.~(\ref{eq:31}) and (\ref{eq:32}) simplify to \cite[Corrolary 17, Corollary 20; p.\ 8]{JM2}, respectively.

\section{Proofs of Four Lemmas}\label{sec:4}

We begin with the proof of Lemma \ref{l:2}. We shall use well-known \cite[The Central Binomial Coefficient, p.\ 26]{Koshy} recurrence relation for the central binomial coefficient:
\begin{equation}
\binom{2(l+1)}{l+1}=\frac{2(2l+1)}{l+1}\binom{2l}{l}.\label{eq:33}
\end{equation}
\subsection{Proof of Lemma \ref{l:2}}

We start from Eq.~(\ref{eq:29}).

By Eq.~(\ref{eq:26}), Eq.~(\ref{eq:29}) becomes
\begin{align*}
R_t(2n,l)&=(n+l+1)\sum_{k=t}^{2n-t}(-1)^k \frac{2l+1}{(k+l+1)(2n-k+l+1)}\binom{2n-2t}{k-t}S(k,l)S(2n-k,l),\\
&=(n+l+1)\sum_{k=t}^{2n-t}(-1)^k (2l+1)\binom{2n-2t}{k-t} \frac{\binom{2k}{k}\binom{2l}{l}}{(k+l+1)\binom{k+l}{k}}\frac{\binom{2(2n-k)}{2n-k}\binom{2l}{l}}{(2n-k+l+1)\binom{2n-k+l}{2n-k}}.
\end{align*}

It is readily verified  \cite[Eq.~(1.4), p.~5]{Koshy}  that 
\begin{align*}
(k+l+1)\binom{k+l}{k}&=(l+1)\binom{k+l+1}{k},\\
(2n-k+l+1)\binom{2n-k+l}{2n-k}&=(l+1)\binom{2n-k+l+1}{2n-k}.
\end{align*}

By using the two equations above, we obtain the equation for the $R_t(2n,l)$ sum:
\begin{equation}
R_t(2n,l)=(n+l+1)\sum_{k=t}^{2n-t}(-1)^k \binom{2n-2t}{k-t} \frac{\binom{2k}{k}(2l+1)\binom{2l}{l}}{(l+1)\binom{k+l+1}{k}}\frac{\binom{2(2n-k)}{2n-k}\binom{2l}{l}}{(l+1)\binom{2n-k+l+1}{2n-k}}.\label{eq:34}
\end{equation}

By using Eq.~(\ref{eq:33}), Eq.~(\ref{eq:34}) becomes
\begin{align}
R_t(2n,l)&=\frac{(n+l+1)}{2}\sum_{k=t}^{2n-t}(-1)^k \binom{2n-2t}{k-t} \frac{\binom{2k}{k}\binom{2(l+1)}{l+1}}{\binom{k+l+1}{k}}\frac{\binom{2(2n-k)}{2n-k}\binom{2l}{l}}{(l+1)\binom{2n-k+l+1}{2n-k}},\notag\\
&=\frac{(n+l+1)}{2}\sum_{k=t}^{2n-t}(-1)^k \binom{2n-2t}{k-t} S(k,l+1)\frac{\binom{2(2n-k)}{2n-k}\binom{2l}{l}}{(l+1)\binom{2n-k+l+1}{2n-k}}.\label{eq:35}
\end{align}
Eq.~(\ref{eq:35}) is equivalent to the following equation:
\begin{equation}
R_t(2n,l)=\frac{(n+l+1)}{4(2l+1)}\sum_{k=t}^{2n-t}(-1)^k \binom{2n-2t}{k-t} S(k,l+1)\frac{\binom{2(2n-k)}{2n-k}2(2l+1)\binom{2l}{l}}{(l+1)\binom{2n-k+l+1}{2n-k}}.\label{eq:36}
\end{equation}

Again, by Eq.~(\ref{eq:33}),  Eq.~(\ref{eq:36}) becomes
\begin{align}
R_t(2n,l)&=\frac{(n+l+1)}{4(2l+1)}\sum_{k=t}^{2n-t}(-1)^k \binom{2n-2t}{k-t} S(k,l+1)\frac{\binom{2(2n-k)}{2n-k}\binom{2(l+1)}{l+1}}{\binom{2n-k+l+1}{2n-k}},\notag\\
&=\frac{(n+l+1)}{4(2l+1)}\sum_{k=t}^{2n-t}(-1)^k \binom{2n-2t}{k-t} S(k,l+1)S(2n-k,l+1).\label{eq:37}
\end{align}

Eqns.~(\ref{eq:19}) and (\ref{eq:37}) complete the proof of Lemma \ref{l:2}.

\subsection{Proof of Lemma \ref{l:3}}

We start from Eq.~(\ref{eq:24}). Note that the last term in  Eq.~(\ref{eq:24}) is equal to zero. Therefore,  Eq.~(\ref{eq:24})  is equivalent to
\begin{equation}
P_t(n,l)=\sum_{k=t}^{n-1-t}(-1)^k(n-t-k)\binom{n-2t}{k-t}S(k,l)S(n-k,l).\label{eq:38}
\end{equation}

It is readily verified  \cite[Eq.~(1.2), p.\ 5]{Koshy} that 
\begin{equation}
(n-t-k)\binom{n-2t}{k-t}=(n-2t)\binom{n-1-2t}{k-t}.\notag
\end{equation}

By the equation above, Eq.~(\ref{eq:38}) becomes
\begin{equation}
P_t(n,l)=(n-2t)\sum_{k=t}^{n-1-t}(-1)^k\binom{n-1-2t}{k-t}S(k,l)S(n-k,l).\label{eq:39}
\end{equation}

Note that, in Eq.~(\ref{eq:39}), $k <n$.

Eq.~(\ref{eq:39}) is equivalent to
\begin{equation}
P_t(n,l)=(n-2t)\sum_{k=t}^{n-1-t}(-1)^k\binom{n-1-2t}{k-t}S(k,l) \frac{\binom{2n-2k}{n-k}\binom{2l}{l}}{\binom{n-k+l}{n-k}}.\label{eq:40}
\end{equation}

It is known that, for positive $n$, $\binom{2n}{n}$ is always an even integer. It is readily verified \cite[Central Binomial Coefficient, p.\ 15]{Koshy} that
\begin{equation}
\binom{2n}{n}=2\binom{2n-1}{n}.\notag
\end{equation}

By the equation above, Eq.~(\ref{eq:40}) becomes
\begin{equation}
P_t(n,l)=2(n-2t)\sum_{k=t}^{n-1-t}(-1)^k\binom{n-1-2t}{k-t}S(k,l) \frac{\binom{2n-2k-1}{n-k}\binom{2l}{l}}{\binom{n-k+l}{n-k}}.\label{eq:41}
\end{equation}

Since $n-k>0$,  Eq.~(\ref{eq:41}) is equivalent to 
\begin{equation}
P_t(n,l)=2(n-2t)\sum_{k=t}^{n-1-t}(-1)^k\binom{n-1-2t}{k-t}S(k,l)\frac{(n-k)\binom{2n-2k-1}{n-k}\binom{2l}{l}}{(n-k)\binom{n-k+l}{n-k}}.\label{eq:42}
\end{equation}

It is readily verified \cite[Eq.~(1.1), p.\ 5]{Koshy} that
\begin{align*}
(n-k)\binom{2n-2k-1}{n-k}&=(2n-2k-1)\binom{2(n-1-k)}{n-1-k},\\
(n-k)\binom{n-k+l}{n-k}&=(n-k+l)\binom{n-1+k+l}{n-1-k}.
\end{align*}

By using the two equations above, Eq.~(\ref{eq:42}) becomes
\begin{equation}
P_t(n,l)=2(n-2t)\sum_{k=t}^{n-1-t}(-1)^k\binom{n-1-2t}{k-t}S(k,l) \frac{2n-2k-1}{n-k+l} \frac{\binom{2(n-1-k)}{n-1-k}\binom{2l}{l}}{\binom{n-1+k+l}{n-1-k}}.\label{eq:43}
\end{equation}

Eq.~(\ref{eq:43}) is equivalent to
\begin{align}
P_t(n,l)&=2(n-2t)\sum_{k=t}^{n-1-t}(-1)^k\binom{n-1-2t}{k-t} \frac{2(n-1-k)+1}{n-k+l} S(k,l) S(n-1-k,l),\notag\\
&=2(n-2t)\sum_{k=t}^{n-1-t}(-1)^k\binom{n-1-2t}{k-t} (2-\frac{2l+1}{n-k+l}) S(k,l) S(n-1-k,l).\label{eq:44}
\end{align}

Note that changing $k$ to $n-k$ in Eq.~(\ref{eq:25}), Eq.~(\ref{eq:25}) becomes 
\begin{equation}
R_t(n,l)=(-1)^n\sum_{k=t}^{n-t}(-1)^k \frac{2l+1}{n-k+l+1}\binom{n-2t}{k-t}S(k,l)S(n-k,l).\label{eq:45}
\end{equation}

By setting $n:=n-1$ in Eq.~(\ref{eq:45}), we obtain that
\begin{equation}
R_t(n-1,l)=(-1)^{n-1}\sum_{k=t}^{n-1-t}(-1)^k \frac{2l+1}{n-k+l}\binom{n-1-2t}{k-t}S(k,l)S(n-1-k,l).\label{eq:46}
\end{equation}

By using Eqns.~ (\ref{eq:19}), and (\ref{eq:46}), Eq.~(\ref{eq:44}) implies an interesting recurrence relation:
\begin{equation}
P_t(n,l)=4(n-2t)\varPsi_t(n-1,l)+(-1)^n2(n-2t)R_t(n-1,l).\label{eq:47}
\end{equation}

By setting $n:=2n$, Eq.~(\ref{eq:47}) becomes
\begin{equation}
P_t(2n,l)=4(2n-2t)\varPsi_t(2n-1,l)+2(2n-2t)R_t(2n-1,l).\label{eq:48}
\end{equation}

By Eq.~(\ref{eq:20}), Eq.~(\ref{eq:48}) becomes
\begin{equation}
P_t(2n,l)=4(n-t)R_t(2n-1,l).\label{eq:49}
\end{equation}

Finally, by Eq.~(\ref{eq:28}), it follows that 
\begin{equation}
(n-t)\varPsi_t(2n,l)=4(n-t)R_t(2n-1,l).\label{eq:50}
\end{equation}

Since $t<n$, by canceling the factor $n-t$ in Eq.~(\ref{eq:50}),  Lemma \ref{l:3} follows, as desired.

\subsection{Proof of Lemma \ref{l:4}}

For the proof of  Lemma \ref{l:4}, we use $T_t(n,l)$ sums.

Eq.~(\ref{eq:27}) is equivalent to
\begin{equation}
T_t(n,l)=(2l+1)\sum_{k=t}^{n-t}(-1)^k \, \bigl{(}\frac{n+l+1-t}{k+l+1}-1\bigr{)}\binom{n-2t}{k-t}S(k,l)S(n-k,l).\notag
\end{equation}

By the equation above, it follows that
\begin{equation}
T_t(n,l)=(n+l+1-t)R_t(n,l)-(2l+1)\varPsi_t(n,l).\label{eq:51}
\end{equation}

By setting $n:=2n-1$ in Eq.~(\ref{eq:51}) and by using  Eq.~(\ref{eq:20}), Eq.~(\ref{eq:51}) becomes
\begin{equation}
T_t(2n-1,l)=(2n+l-t)R_t(2n-1,l)\label{eq:52}.
\end{equation}

On the other side, it can be shown that $T_t(n,l)$ is equal to
\begin{equation}
2(n-2t)(2l+1)\sum_{k=t}^{n-1-t}\frac{2n-2k-1}{(k+l+1)(n-k+l)}\binom{n-1-2t}{k-t}S(k,l)S(n-1-k,l).\label{eq:53}
\end{equation}

The proof of Eq.~(\ref{eq:53}) is similar to the proof of Eq.~(\ref{eq:43}) in the previous subsection. For the sake of brevity, the proof of Eq.~(\ref{eq:53}) is omitted.

The inner sum in  Eq.~(\ref{eq:53}) can be written as
\begin{equation}
\frac{2(n-1-k)+1}{(k+l+1)(n-k+l)}\binom{n-1-2t}{k-t}S(k,l)S(n-1-k,l).\label{eq:54}
\end{equation}

Let 
\begin{equation}
R''_t(n,l)=(2l+1)\sum_{k=t}^{n-t}(-1)^k\frac{n-k}{(k+l+1)(n-k+l+1)}\binom{n-2t}{k-t}S(k,l)S(n-k,l).\label{eq:55}
\end{equation}

By Eqns.~(\ref{eq:55}) and (\ref{eq:54}), Eq.~(\ref{eq:53}) becomes
\begin{equation}
T_t(n,l)=4(n-2t)R''_t(n-1,l)+2(n-2t)R'_t(n-1,l).\label{eq:56}
\end{equation}

By setting $n:=2n-1$ in Eq.~(\ref{eq:56}), it follows that
\begin{equation}
T_t(2n-1,l)=4(2n-1-2t)R''_t(2(n-1),l)+2(2n-1-2t)R'_t(2(n-1),l).\label{eq:57}
\end{equation}

By Eq.~(\ref{eq:55}), changing $k$ to $n-k$, it is readily verified that
\begin{equation}
R''_t(2n,l)=nR'_t(2n,l).\label{eq:58}
\end{equation}

By Eq.~(\ref{eq:58}),  Eq.~(\ref{eq:57}) becomes
\begin{equation}
T_t(2n-1,l)=2(2n-1)(2n-1-2t)R'_t(2(n-1),l).\label{eq:59}
\end{equation}

By Eq.~(\ref{eq:29}), Eq.~(\ref{eq:59}) becomes
\begin{equation}
T_t(2n-1,l)=\frac{2(2n-1)(2n-1-2t)}{n+l}R_t(2n-2,l).\label{eq:60}
\end{equation}

Eqns.~(\ref{eq:52}) and (\ref{eq:60}) complete the proof of Lemma \ref{l:4}.

\subsection{Proof of the Main Lemma}

By Lemmas (\ref{l:2}) and (\ref{l:3}), Eq.~(\ref{eq:32}) of Lemma \ref{l:4} becomes

\begin{equation}
(2n+l-t)\frac{\varPsi_t(2n,l)}{4}=\frac{2(2n-1-2t)(2n-1)}{n+l}\frac{n+l}{4(2l+1)}\varPsi_t(2n-2,l+1)\mbox{.}\label{eq:61}
\end{equation}

Eq.~(\ref{eq:61}) is equivalent to
\begin{align*}
(2n+l-t)\varPsi_t(2n,l)&=\frac{2(2n-1-2t)(2n-1)\varPsi_t(2n-2,l+1)}{2l+1},\\
\frac{2(2n-1-2t)(2n-1)\varPsi_t(2n-2,l+1)}{2l+1}&=(2n+l-t)\varPsi_t(2n,l),\\
\varPsi_t(2n-2,l+1)&=\frac{(2n+l-t)(2l+1)}{2(2n-1-2t)(2n-1)}\varPsi_t(2n,l).
\end{align*}

The last equation above is Eq.~(\ref{eq:23}). This completes the proof of the main lemma \ref{l:1}.

\section{Proof of Theorem \ref{t:2}}\label{sec:5}

Let $n$, $l$, and  $t$ be non-negative integers such that $t\leq n $; and let 
\begin{equation}
\varphi(2n,l,t)=(-1)^t\frac{\binom{2l}{l}\binom{2t}{t}\binom{2(n+l-t)}{n+l-t}\binom{2n}{n}\binom{2n-2t}{n-t}}{\binom{n+l}{n}\binom{2n+l-t}{n}\binom{n}{t}}.\label{eq:62}
\end{equation}

Let us prove that
 \begin{equation}
\varPsi_t(2n,l)=\varphi(2n,l,t).\label{eq:63}
\end{equation}

We shall use induction on $l$ by using the main lemma \ref{l:1}.  By setting $n:=n+1$ in Eq.~(\ref{eq:23}), it follows that
\begin{equation}
\varPsi_t(2n,l+1)=\frac{(2n+2+l-t)(2l+1)}{2(2n+1-2t)(2n+1)}\varPsi_t(2n+2,l).\label{eq:64}
\end{equation}

For the proof of Theorem \ref{t:2}, we shall use Eq.~(\ref{eq:64}), instead of Eq.~(\ref{eq:23}).

By setting $l=0$ in Eq.~(\ref{eq:62}), it follows that
\begin{equation}
\varphi(2n,0,t)=(-1)^t\frac{\binom{2t}{t}\binom{2n}{n}\binom{2n-2t}{n-t}^2}{\binom{2n-t}{n}\binom{n}{t}}.\label{eq:65}
\end{equation}

It is readily verified   \cite[Eq.~(1.4), p.~5]{Koshy}  that 
\begin{equation}
\binom{2n-t}{n}\binom{n}{t}=\binom{2n-t}{t}\binom{2n-2t}{n-t}.\label{eq:66}
\end{equation}

By Eq.~(\ref{eq:66}), Eq.~(\ref{eq:65}) becomes as follows:
\begin{align}
\varphi(2n,0,t)&=(-1)^t\frac{\binom{2t}{t}\binom{2n}{n}\binom{2n-2t}{n-t}^2}{\binom{2n-t}{t}\binom{2n-2t}{n-t}},\notag\\
\varphi(2n,0,t)&=(-1)^t\frac{\binom{2t}{t}\binom{2n}{n}\binom{2n-2t}{n-t}}{\binom{2n-t}{t}}.\label{eq:67}
\end{align}

By Eqns.~(\ref{eq:21}) and (\ref{eq:67}), it follows that Eq.~(\ref{eq:63}) is true for $l=0$. This confirms the induction base.

Furthermore, by setting $t:=n$ in Eq.~(\ref{eq:62}), we obtain that
\begin{align}
\varphi(2n,l,n)&=(-1)^n\frac{\binom{2l}{l}\binom{2n}{n}\binom{2l}{l}\binom{2n}{n}}{\binom{n+l}{n}\binom{n+l}{n}},\notag\\
\varphi(2n,l,n)&=(-1)^n S(n,l)^2.\label{eq:68}
\end{align}

By Eqns.~(\ref{eq:22}) and (\ref{eq:68}), it follows that Eq.~(\ref{eq:63}) is true for $t=n$. For all non-negative integers $n$ and $l$, the equation
\begin{equation}
\varPsi_n(2n,l)=\varphi(2n,l,n)\label{eq:69}
\end{equation}
holds.

Let us prove that the function $\varphi(2n,l,t)$ satisfies recurrence relation (\ref{eq:64}). More precisely, let us prove that
\begin{equation}
\varphi(2n,l+1,t)=\frac{(2n+2+l-t)(2l+1)}{2(2n+1-2t)(2n+1)}\varphi(2n+2,l,t);\label{eq:70}
\end{equation}
where $n$, $l$, and $t$ are non-negative integers such that $t<n$.

By Eq.~(\ref{eq:62}), the right-side of Eq.~(\ref{eq:70}) becomes
\begin{equation}
\frac{(2n+2+l-t)(2l+1)}{2(2n+1-2t)(2n+1)}(-1)^t\frac{\binom{2l}{l}\binom{2t}{t}\binom{2((n+1)+l-t)}{(n+1)+l-t}\binom{2n+2}{n+1}\binom{2n+2-2t}{n+1-t}}{\binom{n+1+l}{n+1}\binom{2n+2+l-t}{n+1}\binom{n+1}{t}}.\label{eq:71}
\end{equation}

By Eq.~(\ref{eq:33}), we know that
\begin{align}
\binom{2n+2}{n+1}=\frac{2(2n+1)}{n+1}\binom{2n}{n},\label{eq:72}\\
\binom{2n+2-2t}{n+1-t}=\frac{2(2n+1-2t)}{n+1-t}\binom{2n-2t}{n-t}.\label{eq:73}
\end{align}

By Eq.~(\ref{eq:72}), Eq.~(\ref{eq:71}) becomes 
\begin{equation}
\frac{(2n+2+l-t)(2l+1)}{2(2n+1-2t)(2n+1)} (-1)^t\frac{\binom{2l}{l}\binom{2t}{t}\binom{2(n+(l+1)-t)}{n+(l+1)-t}\frac{2(2n+1)}{n+1}\binom{2n}{n}\binom{2n+2-2t}{n+1-t}}{\binom{n+1+l}{n+1}\binom{2n+2+l-t}{n+1}\binom{n+1}{t}}.\label{eq:74}
\end{equation}

Eq.~(\ref{eq:74}) simplifies to
\begin{equation}
(-1)^t\frac{(2n+2+l-t)(2l+1)}{2n+1-2t}\frac{\binom{2l}{l}\binom{2t}{t}\binom{2(n+(l+1)-t)}{n+(l+1)-t}\binom{2n}{n}\binom{2n+2-2t}{n+1-t}}{(n+1)\binom{n+1+l}{n+1}\binom{2n+2+l-t}{n+1}\binom{n+1}{t}}.\label{eq:75}
\end{equation}

By Eq.~(\ref{eq:73}), Eq.~(\ref{eq:75}) becomes
\begin{equation}
(-1)^t\frac{(2n+2+l-t)(2l+1)}{2n+1-2t}\frac{\binom{2l}{l}\binom{2t}{t}\binom{2(n+(l+1)-t)}{n+(l+1)-t}\binom{2n}{n}\frac{2(2n+1-2t)}{n+1-t}\binom{2n-2t}{n-t}}{(n+1)\binom{n+1+l}{n+1}\binom{2n+2+l-t}{n+1}\binom{n+1}{t}}.\label{eq:76}
\end{equation}

Eq.~(\ref{eq:76}) is equivalent to
\begin{equation}
(-1)^t (2n+2+l-t)\frac{2(2l+1)\binom{2l}{l}\binom{2t}{t}\binom{2(n+(l+1)-t)}{n+(l+1)-t}\binom{2n}{n}\binom{2n-2t}{n-t}}{(n+1)\binom{2n+2+l-t}{n+1}\binom{n+1+l}{n+1}(n+1-t)\binom{n+1}{t}}.\label{eq:77}
\end{equation}

It is readily verified \cite[Eq.~(1.1), Eq.~(1.2), p.\ 5]{Koshy}  that
\begin{align}
(n+1)\binom{2n+2+l-t}{n+1}&=(2n+2+l-t)\binom{2n+(l+1)}{n},\label{eq:78}\\
(n+1-t)\binom{n+1}{t}&=(n+1)\binom{n}{t}.\label{eq:79}
\end{align}

By Eqns.~(\ref{eq:78}) and (\ref{eq:79}), Eq.~(\ref{eq:77}) becomes
\begin{align}
(-1)^t (2n+2+l-t)\frac{2(2l+1)\binom{2l}{l}\binom{2t}{t}\binom{2(n+(l+1)-t)}{n+(l+1)-t}\binom{2n}{n}\binom{2n-2t}{n-t}}{(2n+2+l-t)\binom{2n+(l+1)-t}{n}\binom{n+1+l}{n+1}(n+1)\binom{n}{t}}\mbox{, or}\notag\\
(-1)^t\frac{2(2l+1)\binom{2l}{l}\binom{2t}{t}\binom{2(n+(l+1)-t)}{n+(l+1)-t}\binom{2n}{n}\binom{2n-2t}{n-t}}{\binom{2n+(l+1)-t}{n}\binom{n+1+l}{n+1}(n+1)\binom{n}{t}}.\label{eq:80}
\end{align}

It is readily verified \cite[Eq.~(1.5), p.\ 5]{Koshy} that 
\begin{equation}
\binom{n+1+l}{n+1}(n+1)=(l+1)\binom{n+l+1}{n}.\label{eq:81}
\end{equation}

By Eq.~(\ref{eq:81}), Eq.~(\ref{eq:80}) becomes
\begin{equation}
(-1)^t\frac{\frac{2(2l+1)}{l+1}\binom{2l}{l}\binom{2t}{t}\binom{2(n+(l+1)-t)}{n+(l+1)-t}\binom{2n}{n}\binom{2n-2t}{n-t}}{\binom{2n+(l+1)-t}{n}\binom{n+(l+1)}{n}\binom{n}{t}}.\label{eq:82}
\end{equation}

By Eq.~(\ref{eq:33}), it follows that
\begin{equation}
\frac{2(2l+1)}{l+1}\binom{2l}{l}=\binom{2(l+1)}{l+1}.\notag
\end{equation}
By the last equation above, Eq.~(\ref{eq:82}) becomes
\begin{equation}
(-1)^t\frac{\binom{2(l+1)}{l+1}\binom{2t}{t}\binom{2(n+(l+1)-t)}{n+(l+1)-t}\binom{2n}{n}\binom{2n-2t}{n-t}}{\binom{2n+(l+1)-t}{n}\binom{n+(l+1)}{n}\binom{n}{t}}.\label{eq:83}
\end{equation}

By Eq.~(\ref{eq:62}), it follows that Eq.~(\ref{eq:83}) is equal to $\varphi(2n,l+1,t)$.

Therefore, Eqns.~(\ref{eq:71}), (\ref{eq:74}), (\ref{eq:75}), (\ref{eq:76}), (\ref{eq:77}), (\ref{eq:80}), (\ref{eq:82}), and (\ref{eq:83}) complete the proof of Eq.~(\ref{eq:70}).

Let $l_0$ be a fixed non-negative integer. Let us assume that the following equation is true for all non-negative integers $n$ and $t$ such that $t\leq n$:
 \begin{equation}
\varPsi_t(2n,l_0)=\varphi(2n,l_0,t).\label{eq:84}
\end{equation}

Eq.~(\ref{eq:84}) is our induction hypothesis.

Let us prove that Eq.~(\ref{eq:84}) implies the following equation:
 \begin{equation}
\varPsi_t(2n,l_0+1)=\varphi(2n,l_0+1,t);\label{eq:85}
\end{equation}
where $n$ and $t$ are non-negative integers such that $t\leq n$.

We treat two cases: $t=n$ and $t<n$ separately.

By setting $l=l_0+1$ in Eq.~(\ref{eq:69}), it follows that Eq.~(\ref{eq:85}) is true for $t=n$.
This proves the first case.

Let us assume that $t<n$. 

By setting $l:=l_0$ in Eq.~(\ref{eq:64}), we obtain that
\begin{equation}
 \varPsi_t(2n,l_0+1)=\frac{(2n+2+l_0-t)(2l_0+1)}{ 2(2n+1-2t)(2n+1)} \varPsi_t(2n+2,l_0).\label{eq:86}
\end{equation}

By setting $l:=l_0$ in Eq.~(\ref{eq:70}), we obtain that
\begin{equation}
\varphi(2n,l_0+1,t)=\frac{(2n+2+l_0-t)(2l_0+1)}{2(2n+1-2t)(2n+1)}\varphi(2n+2,l_0,t).\label{eq:87}
\end{equation}

By setting $n:=n+1$ in Eq.~(\ref{eq:84}) (the induction hypothesis),  it follows that
\begin{equation}
\varPsi_t(2n+2,l_0)=\varphi(2n+2,l_0,t).\label{eq:88}
\end{equation}

 For $t<n$,  by Eqns.~(\ref{eq:86}), (\ref{eq:87}), and (\ref{eq:88}), Eq.~(\ref{eq:85}) follows. 

This proves the second case.

Therefore, we conclude that Eq.~(\ref{eq:84}) implies  Eq.~(\ref{eq:85}). This proves the induction step.
By induction, Theorem \ref{t:2} follows, as desired.
Also this completes the proof of Eq.~(\ref{eq:18}).

\section{Proof of Theorem \ref{t:1}}\label{sec:6}

By setting $t=0$ in Theorem \ref{t:2}, we obtain that
\begin{align}
\sum_{k=0}^{2n}(-1)^k\binom{2n}{k}S(k,l)S(2n-k,l)&=\frac{\binom{2l}{l}\binom{2(n+l)}{n+l}\binom{2n}{n}^2}{\binom{n+l}{n}\binom{2n+l}{n}},\label{eq:89}\\
&=\frac{\binom{2n}{n}\binom{2l}{l}}{\binom{n+l}{n}}\frac{\binom{2(n+l)}{n+l}\binom{2n}{n}}{\binom{2n+l}{n}},\notag\\
&=S(n,l)S(n+l,n).\notag
\end{align}

The last equation above proves Theorem \ref{t:1}.

\begin{remark}\label{r:1}
Theorem \ref{t:1} can be written as
\begin{equation}
\sum_{k=0}^{2n}(-1)^k\binom{2n}{k}S(k,l)S(2n-k,l)
=\frac{\binom{2l}{l}\binom{2(n+l)}{n+l}\binom{2n}{n}}{\binom{2n+l}{l}}.\label{eq:90}
\end{equation}

By symmetry of binomial coefficients and \cite[Eq.~(1.4), p.~5]{Koshy}, it follows that

 \begin{equation}
\binom{2n+l}{n}\binom{n+l}{n}=\binom{2n+l}{l}\binom{2n}{n}.\label{eq:91}
\end{equation}

By Eqns.~(\ref{eq:89}) and (\ref{eq:91}), Eq.~(\ref{eq:90}) follows, as desired.
\end{remark}

\section{Proof of Theorem \ref{t:3}}\label{sec:7}

We consider the following sum:
\begin{equation}
\varPsi(2n,m,l)=\sum_{k=0}^{2n}(-1)^k\binom{2n}{k}^m S(k,l)S(2n-k,l).\label{eq:92}
\end{equation}

By Theorem \ref{t:1} and Eq.~(\ref{eq:92}), we obtain that $\varPsi(2n,1,l)=S(n,l)S(n+l,n)$. Hence Theorem \ref{t:3} is true for $m=1$.

Our proof of Theorem \ref{t:3} consists from two parts.

In the first part, we show that Theorem \ref{t:3} is true for $m=2$. In the second part we prove that $D_{\varPsi}(2n,j,1;l)$ is divisible by $S(n,l)$ for all non-negative integers $n$, $j$, and $l$ such that $j\leq n$. Then, by Eq.~(\ref{eq:12}), it follows that Theorem \ref{t:3} is true for all integers m greater than $2$.

\subsection{The First Part}

By setting $A:=\varPsi$ and $n:=2n$ in Relation (\ref{eq:17}), we obtain that
\begin{equation}
 D_{\varPsi}(2n,j,0;l)=\sum_{u=0}^{n-j}\binom{2n-j}{u}\binom{2n-j-u}{j+u}\varPsi_{j+u}(2n,l).\label{eq:93}
\end{equation}

Let $t$ be a non-negative integer such that $t\leq n$ . Let us prove the following equation 
\begin{equation}
\binom{2n-t}{t}\varPsi_{t}(2n,l)=(-1)^t\frac{\binom{2l}{l}\binom{2t}{t}\binom{2(n+l-t)}{n+l-t}\binom{2n}{n}\binom{2n-t}{n}}{\binom{n+l}{n}\binom{2n+l-t}{n}}.\label{eq:94}
\end{equation}

By Eq.~(\ref{eq:18}), it follows that
\begin{equation}
\binom{2n-t}{t}\varPsi_{t}(2n,l)=\binom{2n-t}{t}(-1)^t\frac{\binom{2l}{l}\binom{2t}{t}\binom{2(n+l-t)}{n+l-t}\binom{2n}{n}\binom{2n-2t}{n-t}}{\binom{n+l}{n}\binom{2n+l-t}{n}\binom{n}{t}}.\label{eq:95}
\end{equation}

Eq.~(\ref{eq:95}) is equivalent to 
\begin{equation}
\binom{2n-t}{t}\varPsi_{t}(2n,l)=(-1)^t \frac{\binom{2l}{l}\binom{2t}{t}\binom{2(n+l-t)}{n+l-t}\binom{2n}{n}}{\binom{n+l}{n}\binom{2n+l-t}{n}} \frac{\binom{2n-t}{t}\binom{2n-2t}{n-t}}{\binom{n}{t}}.\label{eq:96}
\end{equation}

By symmetry of binomial coefficients and \cite[Eq.~(1.4), p.~5]{Koshy}, it follows that  
\begin{equation}
 \frac{\binom{2n-t}{t}\binom{2n-2t}{n-t}}{\binom{n}{t}}=\binom{2n-t}{n}.\label{eq:97}
\end{equation}

 By Eqns.~(\ref{eq:96}) and (\ref{eq:97}), Eq.~(\ref{eq:94}) follows.

By setting $t:=j+u$ in Eq.~(\ref{eq:94}), we obtain that
\begin{equation}
\binom{2n-j-u}{j+u}\varPsi_{j+u}(2n,l)=(-1)^{j+u}\frac{\binom{2l}{l}\binom{2(j+u)}{j+u}\binom{2(n+l-j-u)}{n+l-j-u}\binom{2n}{n}\binom{2n-j-u}{n}}{\binom{n+l}{n}\binom{2n+l-j-u}{n}}.\label{eq:98}
\end{equation}

By Eq.~(\ref{eq:98}), Eq.~(\ref{eq:93}) becomes as follows
\begin{align}
 D_{\varPsi}(2n,j,0;l)&=\sum_{u=0}^{n-j}\binom{2n-j}{u}(-1)^{j+u}\frac{\binom{2l}{l}\binom{2(j+u)}{j+u}\binom{2(n+l-j-u)}{n+l-j-u}\binom{2n}{n}\binom{2n-j-u}{n}}{\binom{n+l}{n}\binom{2n+l-j-u}{n}},\notag\\
&=\frac{\binom{2n}{n}\binom{2l}{l}}{\binom{n+l}{n}}\sum_{u=0}^{n-j}(-1)^{j+u}\frac{\binom{2(j+u)}{j+u}\binom{2(n+l-j-u)}{n+l-j-u}\binom{2n-j}{u}\binom{2n-j-u}{n}}{\binom{2n+l-j-u}{n}},\notag\\
&=S(n,l) \sum_{u=0}^{n-j}(-1)^{j+u}\frac{\binom{2(j+u)}{j+u}\binom{2(n+l-j-u)}{n+l-j-u}\binom{2n-j}{u}\binom{2n-j-u}{n}}{\binom{2n+l-j-u}{n}}.\label{eq:99}
\end{align}

By symmetry of binomial coefficients and \cite[Eq.~(1.4), p.~5]{Koshy}, it follows that
 \begin{equation}
\binom{2n-j}{u}\binom{2n-j-u}{n}=\binom{2n-j}{n}\binom{n-j}{u}.\label{eq:100}
\end{equation}

By Eq.~(\ref{eq:100}), Eq.~(\ref{eq:99}) is equal to
\begin{equation}
 D_{\varPsi}(2n,j,0;l)=S(n,l) \sum_{u=0}^{n-j}(-1)^{j+u}\frac{\binom{2(j+u)}{j+u}\binom{2(n+l-j-u)}{n+l-j-u}\binom{2n-j}{n}\binom{n-j}{u}}{\binom{2n+l-j-u}{n}}.\label{eq:101}
\end{equation}

By setting $j:=0$ in  Eq.~(\ref{eq:101}), we have that
\begin{align}
 D_{\varPsi}(2n,0,0;l)=S(n,l) \sum_{u=0}^{n}(-1)^{u}\frac{\binom{2u}{u}\binom{2(n+l-u)}{n+l-u}\binom{2n}{n}\binom{n}{u}}{\binom{2n+l-u}{n}}\mbox{, or}\notag\\
 D_{\varPsi}(2n,0,0;l)=S(n,l) \sum_{u=0}^{n}(-1)^{u}\binom{2u}{u}S(n,n+l-u)\binom{n}{u}.\label{eq:102}
\end{align}

By Eq.~(\ref{eq:12}), it follows that
 
\begin{equation}
\varPsi(2n,2,l)= D_{\varPsi}(2n,0,0;l).\label{eq:103}
\end{equation}

Hence, by Eqns.~(\ref{eq:102}) and (\ref{eq:103}), it follows that
\begin{equation}
\varPsi(2n,2,l)=S(n,l) \sum_{u=0}^{n}(-1)^{u}\binom{2u}{u}S(n,n+l-u)\binom{n}{u}.\label{eq:104}
\end{equation}

By Eq.~(\ref{eq:104}), it follows that Theorem \ref{t:3} is true for $m=2$.

\begin{remark}\label{r:2}
If $l$ is a positive integer then $\varPsi(2n,2,l)$  is divisible by $2S(n,l)$.

Since $l$ is a positive integer and $u\leq n$ in Eq.~(\ref{eq:104}) , it follows that $n-u+l$ is a positive integer also. Then the integer $S(n,n+l-u)$ must be even. By Eq.~(\ref{eq:104}), it follows that 
$\varPsi(2n,2,l)$  is divisible by $2S(n,l)$.
\end{remark}

\subsection{The Second Part}

Let us  calculate $D_{\varPsi}(2n,j,1;l)$; where $j$ is a non-negative integer such that $j\leq n$.
We shall use Eq.~(\ref{eq:13}).

By  setting $t:=0$, $n:=2n$, and $A:=\varPsi$ in Eq.~(\ref{eq:13}), we obtain that
\begin{equation}
 D_{\varPsi}(2n,j,1;l)=\sum_{u=0}^{n-j}\binom{2n}{j+u}\binom{2n-j}{u}D_{\varPsi}(2n,j+u,0;l).\label{eq:105}
\end{equation}

We use Eq.~(\ref{eq:101}) from the previous subsection.
Eq.~(\ref{eq:101}) is equivalent to
\begin{equation}
 D_{\varPsi}(2n,j,0;l)=(-1)^j S(n,l) \binom{2n-j}{n} \sum_{v=0}^{n-j}(-1)^v\frac{\binom{2(j+v)}{j+v}\binom{2(n+l-j-v)}{n+l-j-v}\binom{n-j}{v}}{\binom{2n+l-j-v}{n}}.\label{eq:106}
\end{equation}

By setting $j:=j+u$ in   Eq.~(\ref{eq:106}), it follows that $ D_{\varPsi}(2n,j+u,0;l)$ is equal to
\begin{equation}
(-1)^{j+u} S(n,l) \binom{2n-j-u}{n} \sum_{v=0}^{n-j-u}(-1)^v\frac{\binom{2(j+u+v)}{j+u+v}\binom{2(n+l-j-u-v)}{n+l-j-u-v}\binom{n-j-u}{v}}{\binom{2n+l-j-u-v}{n}}.\label{eq:107}
\end{equation}

Let $Q(n,j+u,l)$ denote the sum
 \begin{equation}
\sum_{v=0}^{n-j-u}(-1)^v\frac{\binom{2(j+u+v)}{j+u+v}\binom{2(n+l-j-u-v)}{n+l-j-u-v}\binom{n-j-u}{v}}{\binom{2n+l-j-u-v}{n}}.\label{eq:108}
\end{equation}

By Eq.~(\ref{eq:108}), Eq.~(\ref{eq:107}) becomes
\begin{equation}
 D_{\varPsi}(2n,j+u,0;l)=(-1)^{j+u} S(n,l) \binom{2n-j-u}{n}Q(n,j+u,l).\label{eq:109}
\end{equation}

By Eqns.~(\ref{eq:105}) and  Eq.~(\ref{eq:109}),  we obtain that $D_{\varPsi}(2n,j,1;l)$ is equal to
\begin{equation}
(-1)^j S(n,l)\sum_{u=0}^{n-j}(-1)^u \binom{2n-j}{u}\binom{2n}{j+u}\binom{2n-j-u}{n} Q(n,j+u,l).\label{eq:110}
\end{equation}

By  symmetry of binomial coefficients and  \cite[Eq.~(1.4), p.~5]{Koshy}, it follows that
\begin{equation}
\binom{2n}{j+u}\binom{2n-j-u}{n}=\binom{2n}{n}\binom{n}{j+u}.\label{eq:111}
\end{equation}

By Eqns.~(\ref{eq:110}) and (\ref{eq:111}), it follows that
\begin{equation}
D_{\varPsi}(2n,j,1;l)=(-1)^j S(n,l)\sum_{u=0}^{n-j}(-1)^u \binom{2n-j}{u}\binom{n}{j+u}\binom{2n}{n}Q(n,j+u,l).\label{eq:112}
\end{equation}

Note that the number $\binom{2n}{n}Q(n,j+u,l)$ in Eq.~(\ref{eq:112})  is an integer.

By Eq.~(\ref{eq:108}), it follows that
\begin{equation}
\binom{2n}{n}Q(n,j+u,l)=\sum_{v=0}^{n-j-u}(-1)^v\binom{2(j+u+v)}{j+u+v}\binom{n-j-u}{v}\frac{\binom{2n}{n}\binom{2(n+l-j-u-v)}{n+l-j-u-v}}{\binom{2n+l-j-u-v}{n}}.\label{eq:113}
\end{equation}
Namely, by Eq.~(\ref{eq:113}), it follows that  $\binom{2n}{n}Q(n,j+u,l)$ is equal to
\begin{equation}
\sum_{v=0}^{n-j-u}(-1)^v\binom{2(j+u+v)}{j+u+v}\binom{n-j-u}{v}S(n,n+l-j-u-v).\label{eq:114}
\end{equation}

By Eq.~(\ref{eq:114}), the number $\binom{2n}{n}Q(n,j+u,l)$ is always an integer.
By Eq.~(\ref{eq:112}), it follows that $D_{\varPsi}(2n,j,1;l)$ is divisible by $S(n,l)$ for all non-negative integers $n$, $j$, and $l$ such that $j\leq n$.
Then, by Eq.~(\ref{eq:13}) and induction, it can be shown that $D_{\varPsi}(2n,j,t;l)$ is divisible by $S(n,l)$ for all positive integers $t$ and for all non-negative integers $n$, $j$, and $l$ such that $j\leq n$. See \cite[How Does This Method Work, p.\ 7]{JM1}.

By setting $m:=t+2$ in Eq.~(\ref{eq:12}), we obtain that
\begin{equation}
\varPsi(2n,t+2,l)= D_{\varPsi}(2n,0,t;l)\mbox{.}\label{eq:115}
\end{equation}

By Eq.~(\ref{eq:115}), it follows that $\varPsi(2n,t+2,l)$ is divisible by $S(n,l)$ for all positive integers $t$  and for all non-negative integers $n$, $j$, and $l$ such that $j\leq n$.
Since $t\geq 1$, it follows $t+2\geq 3$. Therefore, Theorem \ref{t:3} is true for all non-negative integers $n$ and for all positive integers $m$ such that $m\geq3$. This completes the second part and proves Theorem \ref{t:3}.

\begin{remark}\label{r:3}

Let $l$ be a positive integer. Then  $S(n,n+l-j-u-v)$ is always an even integer.

By Eq.~(\ref{eq:114}), it follows that
number  $\binom{2n}{n}Q(n,j+u,l)$, is an even integer.  By Eq~(\ref{eq:112}), it  follows that  $D_{\varPsi}(2n,j,1;l)$ is divisible by $2S(n,l)$ for all non-negative integers $n$, $j$, and $l$ such that $j\leq n$.
By the method of $D$  sums, it can be shown that $\varPsi(2n,m,l)$  is divisible by $2S(n,l)$ for $m\geq 3$.
Finally, by Theorem \ref{t:1} and Remark \ref{r:2}, it follows that $\varPsi(2n,m,l)$  is divisible by $2S(n,l)$ for all non-negative integers $n$ and for all positive integers $m$ and $l$.
\end{remark}

\begin{remark}\label{r:4}
Let $n$ be a non-negative integer and let $m$ be a positive integer. Then the sum  $\varPsi(2n,m,l)$ is divisible   \cite[Thm.\ 1, Corollary 4; p.\ 2]{JM2} by $\binom{2n}{n}$ for $l=0$ and $l=1$. This is not true for all positive integers $l$. For example,
take $l=2$, $n=4$, and $m=1$.
\end{remark}

\section{Acknowledgments}

I would like to thank my former teacher Aleksandar Stankov-Leko.  Thanks to Vanja Vuji\'{c}
for proofreading this paper.  Also I would like to thank  Professor Jeffrey  Shallit for valuable comments
which helped to improve the article.

\bigskip
\hrule
\bigskip

\noindent 2010 {\it Mathematics Subject Classification}:
Primary  05A10 ; Secondary 05A19.

\noindent\emph{Keywords:} super Catalan number, Catalan number, method of $D$ sums, recurrence relation, induction, alternating binomial sum.

\bigskip
\hrule 
\bigskip

\vspace*{+.1in}

\end{document}